\documentclass[12pt]{article}\textwidth 160mm\textheight 235mm
\usepackage{amsfonts}
\usepackage{amssymb}
\usepackage{graphicx}
\usepackage{color}
\usepackage{tikz}

\def\ve{\varepsilon}

\def\hat{\widehat}
\def\tilde{\widetilde}

\def \dist {{\rm dist}}

\def\epi{{\rm epi\,}}

\def\N{{\cal N}}

\def\O{{\cal O}}

\def\R{I\!\!R}
\def\N{I\!\!N}
\def\B{\mathbb B}

\def\ox{\overline{x}}
\def\oy{\overline{y}}
\def\oz{\overline{z}}

\def\disp{\displaystyle}

\def\lim{\mbox{\rm lim}\,}

\def\Limsup{\mathop{{\rm Lim}\,{\rm sup}}}

\def\tto{\;{\lower 1pt \hbox{$\rightarrow$}}\kern -10pt
\hbox{\raise 2pt \hbox{$\rightarrow$}}\;}
\def\Hat{\widehat}
\def\Tilde{\widetilde}
\def\Bar{\overline}
\def\ra{\rangle}
\def\la{\langle}
\def\ve{\varepsilon}
\def\B{I\!\!B}
\def\h{\hfill\Box}
\def\R{\mathbb{R}}
\def\N{I\!\!N}
\def\ox{\bar{x}}
\def\oy{\bar{y}}
\def\oz{\bar{z}}

\def\gph{\mbox{\rm gph}\,}
\def\epi{\mbox{\rm epi}\,}

\def\sup{\mbox{\rm sup}\,}

\def\h{\hfill\triangle}
\def\dn{\downarrow}
\def\O{\Omega}
\def\ph{\varphi}

\def\oR{\Bar{\R}}
\def\lm{\lambda}
\def\gg{\gamma}

\def\al{\alpha}

\def \N{I\!\!N}
\def\th{\theta}

\def\sce{\setcounter{equation}{0}}
\begin{document}
\vspace*{0.7in}
\begin{center}
{\bf IMPLICIT EULER APPROXIMATION AND OPTIMIZATION OF ONE-SIDED LIPSCHITZIAN DIFFERENTIAL INCLUSIONS}\footnote{This research was partly supported by the National Science Foundation under grant DMS-1007132.}\\[2ex]
BORIS S. MORDUKHOVICH and YUAN TIAN\\[2ex]
Department of Mathematics, Wayne State University, Detroit, MI 48202\\[1ex]
Emails: boris@math.wayne.edu, yuantian@wayne.edu
\end{center}
\small{\bf Abstract.} This paper concerns the study of the generalized Bolza problem governed by differential inclusions satisfying the so-called ``relaxed one-sided Lipschitzian" (ROSL) condition with respect to the state variables subject to various types of nonsmooth endpoint constraints. We construct discrete approximations of differential inclusions with ROSL right-hand sides by using the implicit Euler scheme for approximating time derivatives, and then we justify an appropriate well-posedness of such approximations. Our principal result establishes the uniform approximation of strong local minimizers for the continuous-time Bolza problem by optimal solutions to the implicitly discretized finite-difference systems in the general ROSL setting and even by the strengthen $W^{1,2}$-norm approximation of this type in the case ``intermediate" (between strong and weak minimizers) local minimizers under additional assumptions. Finally, we derive necessary optimality conditions for the discretized Bolza problems via suitable generalized differential constructions of variational analysis. The obtained results on the well-posedness of discrete approximations and necessary optimality conditions allow us to justify a numerical approach to solve the generalized Bolza problem for one-sided Lipschitzian differential inclusions by using discrete approximations constructed via the implicit Euler scheme.\\[1ex]

\newtheorem{Theorem}{Theorem}[section]
\newtheorem{Proposition}[Theorem]{Proposition}
\newtheorem{Remark}[Theorem]{Remark}
\newtheorem{Lemma}[Theorem]{Lemma}
\newtheorem{Corollary}[Theorem]{Corollary}
\newtheorem{Definition}[Theorem]{Definition}
\newtheorem{Example}[Theorem]{Example}
\renewcommand{\theequation}{{\thesection}.\arabic{equation}}
\renewcommand{\thefootnote}{\fnsymbol{footnote}}

\normalsize
\section{Introduction and Problem Formulation}\sce

This paper addresses the following optimization problem $(P)$ of the {\em generalized Bolza type} for dynamic systems governed by constrained differential inclusions:
\begin{equation}\label{cost}
\mbox{minimize}\;J[x]:=\varphi_0\big(x(T)\big)+\int_{0}^{T}f\big(x(t),\dot{x}(t),t\big)dt
\end{equation}
over absolutely continuous trajectories $x:[0,T]\to\mathbb{R}^n$ satisfying the differential inclusion
\begin{equation}\label{di}
\dot{x}(t)\in F\big(x(t),t\big)~\;\mbox{a.e.}~\;t\in[0,T]\;\mbox{ with }\;x(0)=x_0\in\R^n
\end{equation}
subject to the geometric and functional endpoint constraints given by, respectively,
\begin{equation}\label{edc}
x(T)\in\Omega\subset\mathbb{R}^n,
\end{equation}
\begin{equation}\label{ine}
\varphi_{i}(x(T))\le 0\;\mbox{ for }\;i=1,\ldots,m,
\end{equation}
\begin{equation}\label{equ}
\varphi_{i}(x(T))=0\;\mbox{ for }\;i=m+1,\ldots,m+r.
\end{equation}
Here $x_0$ is a fixed $n$-vector, $F\colon\R^n\times[0,T]\tto\R^n$ is a set-valued mapping/multifunction, $\O$ is an nonempty set, $f$ and $\ph_i$ for $i=0,\ldots,m+r$ are real-valued functions. Differential inclusion problems of type ($P$) have been well recognized in dynamic optimization and control theory as a convenient framework to cover the vast majority of conventional and nonconventional models arising in optimization and control of dynamical systems described via time derivatives. We refer the reader to the books \cite{bs2,v00} and the bibliographies therein for more discussions, historical overviews as well as applied models governed by differential inclusions. In particular, the differential inclusion model in (\ref{di}) encompasses ODE control systems represented in the parameterized control form
\begin{equation}\label{cont}
\dot x(t)=g(x(t),u(t),t)\;\mbox{ a.e. on }\;[0,T]\;\mbox{ with }\;u(t)\in U(x(t),t),
\end{equation}
where $g\colon\R^n\times U\times[0,T]\to\R^n$ and $U(x,t)\subset U$ is a variable control set belonging to the given control space $U$. Indeed, we can take
\begin{eqnarray}\label{cont-rep}
F(x,t):=\{v\in\R^n|\;v=g(x,u,t)\;\mbox{ for some }\;u\in U(x,t)\}
\end{eqnarray}
in (\ref{di}) to describe (\ref{cont}), omitting here details regarding measurable selections. Observe that the differential inclusion framework (\ref{di}) with the velocity mapping $F$ in (\ref{cont-rep}) covers not only standard control systems with constant control sets but also significantly more challenging problems with {\em feedback} reflected by the dependence of the control sets in (\ref{cont}) on state variables. On the other hand, the optimization problem $(P)$ is {\em intrinsically nonsmooth}, (even when $\O=\R^n$ and all the functions in (\ref{cost}), (\ref{ine}), (\ref{equ}), and (\ref{cont}) are smooth) due to set-valued dynamic constraints in (\ref{di}). Thus the usage and development of appropriate tools of variational analysis and generalized differentiation  are required for the study and applications of $(P)$ and related problems governed by differential inclusions.\vspace*{0.05in}

The {\em method of discrete approximations} has been well recognized as an efficient approach to investigate differential inclusions and optimization problems for them from both qualitative and quantitative/numerical viewpoints; see, e.g., the surveys and books \cite{dl,lv,bs2,s} and the bibliographies therein. A principal question arising in all the aspects and modifications of this method (even without applications to optimization) is about the possibility to approximate, in a suitable sense, feasible trajectories of the given differential inclusion by those for finite-difference inclusions that appear by using one or another scheme to replace time derivatives. The majority of the results in this direction concern explicit Euler schemes under the Lipschitz continuity of velocity mappings with respect to state variables; see \cite{a,dl,lv,bs95,bs2,s} for more details and references.

The other lines of research on discrete approximations of differential inclusions via the explicit Euler schemes invoke the replacement of the Lipschitz continuity of velocity mappings by various {\em one-sided Lipschitzian} conditions; see, e.g., \cite{td91,td02,tdf98,tdeb,lv}. Conditions of this type essentially weaken, from one side, the classical Lipschitz continuity, while from the other side they encompass {\em dissipativity} properties widely used in nonlinear analysis and the theory of monotone operators. Note to this end the so-called {\em modified one-sided Lipschitzian} (MOSL) condition introduced and applied in \cite{tdeb} to justify a certain strong approximation of solution sets for differential inclusions by finite-difference ones obtained via the explicit Euler scheme and to derive in this way a Bogolyubov-type density theorem for the Bolza problem $(P)$ and the corresponding convergence of discrete optimal solutions.\vspace*{0.05in}

In this paper we exploit a weaker property than MOSL known as the {\em relaxed one-sided Lipschitz} (ROSL) condition; see below. The ROSL property of set-valued mappings was introduced by Tzanko Donchev in \cite{td91} under a different name and has already been employed in the studies of various aspects of analysis of set-valued mappings, differential inclusions, and their discrete approximations; see, e.g., \cite{br,td02,tdf98,dfr}. In particular, the paper \cite{tdf98} contains an extension to the ROSL case of the fundamental Filippov theorem about relationships between trajectories and ``quasitrajectories" of Lipschitzian differential inclusions and provides applications of this result to the stability analysis of the explicit Euler scheme. In \cite{br}, similar and related solvability and stability results were developed for the parameterized {\em implicit Euler scheme}
\begin{eqnarray}\label{imp-s}
\Phi_h(x):=\{y\in\R^n|\;y\in x+h F(y)\},\quad h>0,
\end{eqnarray}
generated by ROSL mappings $F$ with compact and convex values. Note that the implicit framework of (\ref{imp-s}) is essentially more involved in comparison with the explicit one
\begin{eqnarray}\label{exp-e}
\Phi_h(x):=\{y\in\R^n|\;y\in x+h F(x)\},\quad h>0,
\end{eqnarray}
studied and applied in \cite{td02,tdf98,dfr} and other publications.

The main goal of this paper is to use the implicit Euler scheme (\ref{imp-s}) to construct and investigate the {\em discrete approximations}
\begin{eqnarray}\label{imp-d}
x^{k}_{j+1}\in x^{k}_{j}+h_{k}F(x^{k}_{j+1},t_{j+1}),\quad k\in\N:=\{1,2,\ldots\}\;\mbox{ with }\;h_k\dn 0\;\mbox{ as }\;k\to\infty,
\end{eqnarray}
of the ROSL differential inclusion (\ref{di}) and the generalized Bolza problem $(P)$ for it with establishing the strong convergence of discrete approximations (in the sense specified below) and deriving necessary optimality conditions for their optimal solutions. To the best of our knowledge, the results obtained in what follows are new for discrete approximations constructed via the implicit Euler scheme even for the case of unconstrained differential inclusions satisfying the classical Lipschitz condition with respect to state variables.\vspace*{0.05in}

After recalling the basic definitions and some background material in Section~2, we develop new results in the aforementioned directions outlined in what follows.

Section~3 presents a constructive procedure allowing us to {\em strongly approximate} in the norm topology of $C[0,T]$ a given feasible trajectory $\ox(\cdot)$ of the differential inclusion (\ref{di}) by feasible solutions to the {\em implicit Euler} finite-difference inclusions (\ref{imp-d}) piecewise linearly extended to $[0,T]$. Furthermore, we justify here even {\em stronger} $W^{1,2}[0,T]$-norm approximation of $\ox(\cdot)$ by feasible extended discrete trajectories in the following two major cases: either $F$ is {\em ROSL and locally graph-convex}, or $F$ is {\em locally Lipschitzian}. Some counterparts of this result involving the {\em explicit} Euler scheme (\ref{exp-e}) can be found (with different proofs) in \cite{bs95,bs2} for Lipschitzian differential inclusions and in \cite{tdeb} for those satisfying the MOSL condition. We are not familiar with any results of this type (involving either the $C[0,T]$ or $W^{1,2}[0,T]$ convergence) for discrete approximations of differential inclusions based on the implicit Euler scheme.

In Section~4 we construct a sequence of finite-difference Bolza type problems $(P_k)$ as $k\in\N$ with the dynamic constraints given by the implicit Euler scheme (\ref{imp-d}) under appropriate approximations of the cost functional (\ref{cost}) and the endpoint constraints in (\ref{edc})--(\ref{equ}). Then we show that optimal solutions to $(P_k)$ and their slight modifications exist for all large $k\in\N$ and {\em norm converge} in the $C[0,T]$ topology for the case of {\em strong local minimizers} of $(P)$ in the general ROSL setting and the $W^{1,2}[0,T]$ topology in the case of {\em intermediate local minimizers} of $(P)$ under the additional assumptions on these minimizers imposed in Section~3. The obtained results seem to be the first achievements in this direction for the {\em implicit} Euler scheme (\ref{imp-d}). It is worth mentioning however that our approach to the strong approximation and convergence results obtained in Sections~3 and 4 require, along with the ROSL condition on the differential inclusion, the {\em unform boundedness} of the velocity sets. This does not allow us to cover the corresponding developments presented of  \cite{chhm,chhm1} for discrete approximations of control problems governed by Moreau's sweeping process, which is described by a dissipative while intrinsically unbounded differential inclusion studied in \cite{chhm} via the explicit Euler scheme by exploiting certain specific features of the sweeping process generated by controlled moving sets.

In addition to the well-posedness results for discrete approximations of $(P)$ via the implicit Euler scheme obtained in Sections~3 and 4, we derive in Section~5 under fairly mild assumptions necessary conditions for optimal solutions to the nonsmooth discrete approximations problems $(P_k)$ associated with the implicit discrete inclusions (\ref{imp-d}) that are different from necessary optimality conditions for the corresponding problems associated with explicit Euler counterparts. These conditions are expressed in terms of the advanced tools of generalized differentiation in variational analysis discussed in Section~2. Due to the established convergence of discrete optimal solutions, the necessary optimality conditions for problems $(P_k)$ obtained in this way can be treated as {\em suboptimality} (almost optimality) conditions for the original Bolza problem $(P)$ and can be also viewed as a certain justification of {\em numerical algorithms} based on discrete approximations. The final Section~6 presents concluding remarks on some topics of further research including deriving necessary optimality conditions for the one-sided Lipschitzian generalized Bolza problem $(P)$ by using the method of discrete approximations.

\section{Basic Definitions and Preliminaries}\sce

Throughout the paper we use standard notation and terminology of variational analysis and generalized differentiation; see, e.g., \cite{bs1,rw}. Recall that $\R^n$ denotes the $n$-dimensional space with the Euclidean norm $|\cdot|$ and the closed unit ball $\B$ and that $\mathcal{CC}(\mathbb{R}^n)$ signifies the space of convex and compact subsets of $\mathbb{R}^n$ endowed with the Pompieu-Hausdorff metric. The distance function associated with an nonempty closed set $\O\subset\R^n$ is denoted by
\begin{eqnarray*}
\mbox{dist}(x,\O):=\disp\min_{y\in\O}|x-y|,\quad x\in\R^n,
\end{eqnarray*}
and the distance between two closed sets $\O_1,\O_2\subset\R^n$ is given by
\begin{eqnarray}\label{dist}
\mbox{dist}(\O_1,\O_2):=\max\Big\{\disp\max_{x\in\O_1}\mbox{dist}(x,\O_2),\,\disp\max_{y\in\O_2}\mbox{dist}(y,\O_1)\Big\}.
\end{eqnarray}
Finally, for an arbitrary a set-valued mapping $F\colon\R^n\tto\R^m$, the Painlev\'e-Kuratowski outer limit of $F$ as $x\to\ox$ is defined by
\begin{eqnarray}\label{pk}
\Limsup_{x\to\ox}F(x):=\Big\{y\in\R^m\Big|\;\exists\,x_k\to\ox,\;y_k\to y\;\mbox{ with }\;y_k\in F(x_k),\;k\in\N\Big\}.
\end{eqnarray}

The following property introduced in \cite{td91} is our standing assumption on the right-hand side $F(\cdot,t)$ of the differential inclusion in (\ref{di}) playing a crucial role in this paper.

\begin{Definition}[relaxed one-sided Lipschitzian condition]\label{rosl} A set-valued mapping $F:\mathbb{R}^n\to\mathcal{CC}(\mathbb{R}^n)$ is called to be {\sc relaxed one-sided Lipschitzian (ROSL)} with constant $l\in\mathbb{R}$ if for any given $x_1,x_2\in\mathbb{R}^n$ and $y_1\in F(x_1)$ there exits $y_2\in F(x_2)$ such that
\begin{equation}\label{osl}
\la y_1-y_2,x_1-x_2\ra\le l|x_1-x_2|^2.                                                                                                                                                                                                                         \end{equation}
\end{Definition}

Note that the number/modulus $l$ in (\ref{osl}) is not required to be positive as in the classical Lipschitz continuity. The ROSL condition is dramatically weaker the standard Lipschitz continuity and essentially relaxes dissipativity and other one-sided Lipschitzian properties; see more discussions and examples in \cite{td91,td02,tdf98,lv}.

The next result on the solvability of the implicit Euler scheme (\ref{imp-s}) under the ROSL condition is taken from \cite[Theorem~4]{br} (the proof of which is based on the Kakutani fixed-point theorem) and is useful in what follows. Recall that a set-valued mapping $F$ is upper semicontinuous (usc) on $\R^n$ if for any $\ox\in\R^n$ and any $\varepsilon>0$ there exists $\gg>0$ such that $F(x)\subset F(\ox)+\varepsilon\B$ whenever $|x-\ox|\le\gg$.

\begin{Lemma}[solvability of the implicit Euler scheme]\label{brth} Let $F\colon\R^n\to\mathcal{CC}(\mathbb{R}^n)$ be usc and ROSL on $\R^n$  with constraint $l\in\R$ such that $lh<1$. Then for any $x,y\in\R^n$, there exists a solution $\bar{y}\in\Phi_{h}(x)$ of the implicit Euler scheme {\rm(\ref{imp-s})} satisfying the estimate
$$
|\bar{y}-y|\le\frac{1}{1-l h}{\rm{dist}}(y,x+hF(y)).
$$
\end{Lemma}

Following \cite{bs95}, we say that a feasible solution $\ox(\cdot)$ to $(P)$ is an {\em intermediate local minimizer} (i.l.m.) of rank $p\in [1,\infty)$ for this problem if there are positive numbers $\varepsilon,\alpha$ such that $J[\bar x]\le J[x]$ for any other feasible solutions $x(\cdot)$ to {\rm($P$)} satisfying the conditions
\begin{eqnarray}\label{ilme}
|x(t)-\bar{x}(t)|<\varepsilon\;\mbox{ as }\;t\in[0,T]\quad\;\mbox{and}\;\quad\alpha\int_{0}^{T}|\dot{x}(t)-\dot{\bar{x}}(t)|^{p}dt<\varepsilon.
\end{eqnarray}

The case of $\alpha=0$ in (\ref{ilme}) corresponds to the classical notion of {\em strong} local minimum and surely includes global solutions to $(P)$. The notion of {\em weak} local minimum corresponds to (\ref{ilme}) with $\alpha\ne 0$ and $p=\infty$; see \cite{bs95,bs2} for detailed discussions and examples.

In what follows we need a certain modification of the i.l.m.\ notion formulated above, which related to some local relaxation stability of the initial problem $(P)$. Along with $(P)$, consider its extended/relaxed version constructed in the line well understood in the calculus of variations and optimal control. Let
\begin{equation}\label{ext-f}
f_{F}(x,v,t):=f(x,v,t)+\delta(v,F(x,t)),
\end{equation}
where $\delta(\cdot,\Lambda)$ is the indicator function of the set $\Lambda$ equal to 0 on $\Lambda$ and to $\infty$ otherwise. Denote by $\hat{f}_{F}(x,v,t)$ the convexification for $f_F$ in the $v$ variable, i.e., the largest convex function majorized by $f_{F}(x,\cdot,t)$ for each $x$ and $t$. The {\em relaxed generalized Bolza problem} $(R)$ consists of minimizing the functional
\begin{equation}\label{rel-cost}
\hat{J}[x]:=\varphi(x(T))+\int_{0}^{T}\hat{f}_{F}(x(t),\dot{x}(t),t)dt
\end{equation}
over absolutely continuous trajectories $x\colon[0,T]\to\R^n$ of the differential inclusion (\ref{di}) with the endpoint constraints (\ref{edc})--(\ref{equ}).

Note that, due to our standing assumption on the convexity and compactness of the values $F(x,t)$ in the differential inclusion (\ref{di}), the validity of the dynamic constraint (\ref{di}) is automatic for any absolutely continuous function $x(\cdot)$ with $\hat{J}[x]<\infty$ in (\ref{rel-cost}). Thus the relaxed problem $(R)$ reduces to the original one $(P)$ if the integrand $f$ is convex with respect to the velocity variable $v$; in particular, when $f$ does not depend on $v$. Furthermore, a remarkable fact for the continuous-time problems under consideration consists of the equality between the infimum values of the coast functionals in $(P)$ and $(R)$, without taking endpoint constraints into account, even when $f$ is not convex in $v$. This fact is known as ``hidden convexity" of continuous-valued variational and control problems and relates to the fundamental results of Bogolyubov's and Lyapunov's types; see, e.g., the books \cite{abm,bs2,v00} for exact formulations and more discussions. The most recent extended version of the Bogolyubov theorem for differential inclusion problems of type $(P)$ was obtained in \cite{tdeb} under the MOSL condition on $F(\cdot,t)$ mentioned in Section~1. This discussion makes more natural the following notion taken from \cite{bs95}.

\begin{Definition}[relaxed intermediate local minimizers]\label{rilm} A feasible solution $\bar{x}(\cdot)$ to the original problem $(P)$ is called a {\sc relaxed intermediate local minimum} $(r.i.l.m.)$ of rank $p\in[1,\infty)$ for $(P)$ if it provides an intermediate local minimum of rank $p$ to the relaxed problem $(R)$ and satisfies the condition $J[\bar x]=\hat{J}[\ox]$.
\end{Definition}

Finally in this section, we recall and briefly discuss the generalized differential constructions of variational analysis introduced by the first author \cite{bs76} and employed in Section~5 for deriving necessary optimality conditions; see the books \cite{bs1,rw} for more details and references on these and related constructions. Given a set $\O\subset\R^n$ locally closed around $\ox\in\O$, the {\em normal cone} to $\O$ at $\ox$ is defined via the outer limit (\ref{pk}) by
\begin{eqnarray}\label{mnc}
N(\bar x;\Omega)=\Limsup_{x\to\bar x}\Big[{\rm{cone}}(x-\Pi(x,\Omega))\Big],
\end{eqnarray}
where $\Pi(x,\Omega)=\big\{w\in\Omega\;\mbox{s.t.}\;|x-w|=\dist(x,\Omega)\big\}$ is the Euclidean projector of $x$ on $\O$, and where the symbol ``cone" stands for the conic hull of the set in question. This normal cone reduces to the classical normal cone of convex analysis when $\O$ is convex, which it may take nonconvex values in rather simple situations as, e.g., for the graph of the function $|x|$ and the epigraph of the function $-|x|$ on $\R$. Nevertheless the normal cone (\ref{mnc}) and the related generalized differential constructions for functions and mappings enjoy comprehensive calculus rules based on variational/extremal principles of variational analysis; see \cite{bs1,rw} and the references therein.

Given now an extended-real-valued and lower semicontinuous function $\ph\colon\R^n\to\oR:=(-\infty,\infty]$ finite at $\ox$, we define its {\em subdifferential} at $\ox$ geometrically
\begin{eqnarray}\label{subd}
\partial\ph(\ox):=\Big\{v\in\R^n\Big|\;(v,-1)\in N((\ox,\ph(\ox));\epi\ph)\Big\}
\end{eqnarray}
via the normal cone (\ref{mnc}) to the epigraphical set
$$
\epi\ph:=\{(x,\al)\in\R^n\times\R|\;\al\ge\ph(x)\}
$$
of $\ph$. The reader can find in \cite{bs1,rw} various analytical representations and properties of the subgradient mapping $\partial\ph\colon\R^n\tto\R^n$ used in what follows.

We recall also the {\em symmetric subdifferential} construction for a continuous function $\ph\colon\R^n\to\R$  at $\ox$ defined by
\begin{eqnarray}\label{sym}
\partial^{0}\ph(\ox):=\partial\ph(\ox)\cup(-\partial(-\ph)(\ox))
\end{eqnarray}
and employed in Section~5 for expressing necessary optimality conditions for equality constraints. Note the symmetry relation
$$
\partial^0(-\ph)(\ox)=-\partial^0\ph(\ox),
$$
which does not hold for the unilateral subdifferential construction (\ref{subd}).

\section{Strong Approximation via Implicit Euler Scheme}\sce

In this section we justify the possibility to {\em strongly approximate} (in the norm topology of either $C[0,T]$ or $W^{1,2}[0,T])$ feasible trajectories of the ROSL inclusion (\ref{di}) constructed via the implicit Euler scheme. Given an {\em arbitrary} trajectory $\ox(\cdot)$ of (\ref{di}), we impose the following assumptions of $F$ near $\ox(\cdot)$ standing throughout Sections~3 and 4. For simplicity, suppose that the uniform boundedness and ROSL moduli below are constant on $[0,T]$. They can obviously be replaced by the continuous functions $m_F(t)$ and $l(t)$ on this compact interval while the proofs of the main results presented in Sections~3 and 4 can be modified to more general cases of the Riemann and Lebesgue integrability.

\begin{enumerate}
\item[(\bf H1)]\label{h1}
There exists an open set $U\subset\mathbb{R}^n$ and a number $m_F>0$ such that $\bar{x}(t)\in U$ for all $t\in[0,T]$ and the multifunction $F\colon U\times[0,T]\to{\cal CC}(\R^n)$ from {\rm(\ref{di})} satisfies the {\em uniform boundedness condition}
\begin{eqnarray*}
F(x,t)\subset m_{F}\B\;\mbox{ for all }\;x\in U,\;\mbox{ a.e. }\;t\in[0,T].
\end{eqnarray*}

\item[(\bf H2)]\label{h1u}
Given $U$ from {\rm (H1)}, for all $x_1,x_2\in U$, a.e.\ $t\in[0,T]$, and $y_1\in F(x_1,t)$ there exists $y_2\in F(x_2,t)$ such that we have the {\em relaxed one-sided Lipschitzian condition}
\begin{eqnarray*}
\la y_1-y_2,x_1-x_2\ra\le l|x_1-x_2|^2.
\end{eqnarray*}

\item[(\bf H3)]\label{h2}
The multifunction $F(\cdot,t)$ is {\em continuous} on the neighborhood $U$ from {\rm(H1)} for a.e.\ $t\in[0,T]$ while $F(x,\cdot)$ is {\em a.e.\ continuous} on $[0,T]$ uniformly in $x\in U$ with respect to the Pompieu-Hausdorff metric.
\end{enumerate}

We now construct a finite-difference approximation of the differential inclusion in (\ref{di}) by using the {\em implicit Euler method} to replace the time derivative by
$$
x(t+h)\in x(t)+hF(x(t+h),t)\;\mbox{ as }\;h\dn 0.
$$
To formalize this process, for any $k\in\N$ define the discrete grid/mesh on $[0,T]$ by
$$
T_{k}:=\big(t_{j}\big|\;j=0,1,\ldots,k\big)\;\mbox{ with }\;t_{0}:=0,\;t_{k}:=T,\;\mbox{ and stepsize }\;h_{k}:=T/k=t_{j+1}-t_{j}.
$$
Then the corresponding discrete inclusions associated with (\ref{di}) via the implicit Euler scheme are constructed as follows:
\begin{equation}\label{ie}
x^{k}_{j+1}\in x^{k}_{j}+h_{k}F(x^{k}_{j+1},t_{j+1})\;\mbox{ for }\;j=0,\ldots,k-1,
\end{equation}
where the starting vector $x_0$ in (\ref{ie}) is taken from (\ref{di}).\vspace*{0.05in}

The next theorem justifies the aforementioned strong $W^{1,2}[0,T]$-approximation of feasible solutions to (\ref{di}) by those for the discrete inclusions (\ref{ie}).

\begin{Theorem}[discrete approximation of ROSL differential inclusions]\label{sa} Let $\bar{x}(\cdot)$ be a feasible trajectory for {\rm(\ref{di})} such that $\dot{\ox}(t)$ is Riemann integrable on $[0,T]$  and the standing assumptions {\rm (H1)--(H3)} are satisfied. Then the following assertions hold:

{\bf (i)} There is a sequence $\{z^{k}_{j}|j=0,\ldots,k\}$ of feasible solutions to the discrete inclusions {\rm(\ref{ie})} such that their piecewise linearly extensions to $[0,T]$ converge to $\ox(t)$ uniformly on $[0,T]$, i.e., in the norm topology of $C[0,T]$.

{\bf (ii)} Assume in addition that either the graph of $F(\cdot,t)$ is locally convex around $(\ox(t),\dot{\ox}(t))$, or $F(\cdot,t)$ is locally Lipschitzian around $\ox(t)$ for a.e.\ $t\in[0,T]$. Then there is a sequence $\{z^{k}_{j}|j=0,\ldots,k\}$ of feasible solutions to {\rm(\ref{ie})} such that the piecewise constantly extended to $[0,T]$ discrete velocity functions
\begin{equation}\label{dis-vel}
v^{k}(t):=\frac{z^{k}_{j+1}-z^{k}_{j}}{h_k},\quad t\in(t_{j},t_{j+1}],\quad j=0,\ldots,k-1,
\end{equation}
converge to $\dot{\bar{x}}(\cdot)$ as $k\to\infty$ in the norm topology of $L^{2}[0,T]$, which is equivalent to the $W^{1,2}[0,T]$-norm convergence of the piecewise linear functions $z^k(t)$ represented by
\begin{eqnarray}\label{dis-tr}
z^k(t)=x_0+\disp\int_{0}^t v^k(s)\,ds\;\mbox{ for all }\;t\in[0,T],\quad k=1,2,\ldots.
\end{eqnarray}
\end{Theorem}
{\bf Proof.} Fix an arbitrary feasible trajectory $\ox(t)$ for (\ref{di}) from the formulation of the theorem and denote $\ox_j:=\ox(t_j)$. Taking into account the density of step functions in $L^{1}[0,T]$, we can find without loss of generality a sequence of functions $w^{k}(\cdot)$ on $[0,T]$ such that $w^{k}(t)$ are constant on $(t_{j},t_{j+1}]$ and $w^{k}(t)$ converge to $\dot{\bar{x}}(t)$ as $k\rightarrow\infty$ in the norm topology of $L^{1}[0,T]$. It follows from (H1) that
\begin{eqnarray*}
|w^k(t)|\le m_F+1\;\mbox{ for all }\;t\in[0,T]\;\mbox{ and }\;k\in\N.
\end{eqnarray*}
Define further the sequences
\begin{equation}\label{xi}
w^k_j:=w^k(t_j)\;\mbox{ for }\;j=1,\ldots,k\;\mbox{ and }\;\xi_k:=\int^T_0|\dot{\bar x}(t)-w^k(t)|\,dt\to 0\;\mbox{as}\;k\to\infty
\end{equation}
and for each $k\in\N$ form recurrently the collection of vectors $\{y^k_0,\ldots,y^k_k\}$ by
\begin{equation}\label{y}
y^{k}_{j+1}:=y^{k}_{j}+h_{k}w^{k}_{j+1}\;\mbox{ for }\;j=0,\ldots,k-1\;\mbox{ with }\;y^{k}_{0}=x_0.
\end{equation}
Note that the continuous-time vector functions
$$
y^{k}(t):=x_0+\disp\int^{t}_{0}w^{k}(s)\,ds,\quad 0\le t\le T,
$$
are piecewise linear extensions of the discrete ones (\ref{y}) on $[0,T]$ satisfying the estimate
\begin{eqnarray}\label{y1}
|y^{k}(t)-\bar{x}(t)|\le\int^{t}_{0}|w^{k}(s)-\dot{\bar x}(s)|\,ds\le\xi_{k}\;\mbox{ for all }\;t\in [0,T],\;k\in\N,
\end{eqnarray}
where $\xi_k$ is taken from (\ref{xi}). Now we construct a sequence of discrete trajectories for (\ref{ie}) by the following {\em algorithmic procedure}.

To define such trajectories $z^k=(z^k_0,\ldots,z^k_k)$ of (\ref{ie}), put $z^{k}_{0}:=x_0$ and suppose that the vectors $z^{k}_{j}$ have been already found. Then for any $k\in\N$ sufficiently large (i.e., when $h_k$ is small) we use the solvability result from Lemma~\ref{brth} valid under assumptions (H2) and (H3) and solve the discrete inclusions (\ref{ie}) for $z^{k}_{j+1}$. Taking into account the error estimate in Lemma~\ref{brth}, the construction of $y^k_j$ in (\ref{y}), and the corresponding properties of the distance (\ref{dist}), we deduce that the vector $z^k_{j+1}$ satisfies the discrete inclusion
\begin{equation}\label{iez}
z^{k}_{j+1}\in z^{k}_{j}+h_{k}F(z^{k}_{j+1},t_{j+1})
\end{equation}
and the following relationships for each $j\in\{1,\ldots,k-1\}$ and all (large) $k\in\N$:
\begin{eqnarray*}
|z^{k}_{j+1}-y^{k}_{j+1}|&\le&\frac{1}{1-l h_k}\dist\big(y^{k}_{j+1},z^{k}_{j}+h_{k}F(y^{k}_{j+1},t_{j+1})\big)\\
&\le&\frac{1}{1-l h_k}\dist(y^{k}_{j+1},y^{k}_{j}+h_{k}F(y^{k}_{j+1},t_{j+1}))\\
&+&\frac{1}{1-l h_k}\dist\big(y^{k}_{j}+h_k F(y^{k}_{j+1},t_{j+1}),z^{k}_{j}+h_{k}F(y^{k}_{j+1},t_{j+1})\big)\\
&\le&\frac{|z^{k}_{j}-y^{k}_{j}|}{1-l h_k}+\frac{h_k}{1-l h_k}\dist\Big(\frac{y^{k}_{j+1}-y^{k}_{j}}{h_k},F(y^{k}_{j+1},t_{j+1})\Big)\\
&=&\frac{|z^{k}_{j}-y^{k}_{j}|}{1-l h_k}+\frac{h_k}{1-l h_k}\dist\big(w^{k}_{j+1},F(y^{k}_{j+1},t_{j+1})\big).
\end{eqnarray*}
Thus we arrive at the estimate valid for all $j=0,\ldots,k-1$ and $k\in\N$:
\begin{eqnarray}\label{zjyj}
|z^{k}_{j+1}-y^{k}_{j+1}|\le\frac{|z^{k}_{j}-y^{k}_{j}|}{1-l h_k}+\frac{h_k}{1-l h_k}\dist\big(w^{k}_{j+1},F(y^{k}_{j+1},t_{j+1})\big).
\end{eqnarray}
Proceeding further by induction implies that
$$
|z^{k}_{j+1}-y^{k}_{j+1}|\le h_{k}\sum^{j+1}_{m=1}\Big(\frac{1}{1-lh_k}\Big)^{j+2-m}\dist\big(w^{k}_{m},F(y^{k}_{m},t_{m})\big),
$$
which yields by choosing $k\in\N$ with $lh_{k}<1/2$ that
\begin{eqnarray}\label{zy}
\begin{array}{ll}
|z^{k}_{j+1}-y^{k}_{j+1}|&\le h_{k}\disp\sum^{j+1}_{m=1}(1+2lh_{k})^{j+2-m}\dist\big(w^{k}_{m},F(y^{k}_{m},t_{m})\big)\\
&\le h_{k}e^{2lT}\disp\sum^{j+1}_{m=1}\dist(w^{k}_{m},F(y^{k}_{m},t_{m})).
\end{array}
\end{eqnarray}

We recall next the average modulus of continuity of $F$ defined by
\begin{eqnarray*}
\tau(F;h):=\disp\sup_{x\in U}\int^{T}_{0}\sup\Big\{\dist(F(x,t'),F(x,t''))\Big|\;t',t''\in\Big[t-\frac{h}{2},t+\frac{h}{2}\Big]\Big\}\,dt
\end{eqnarray*}
and consider the quantities $\zeta_k$ with the estimates
\begin{eqnarray}\label{zeta}
\begin{array}{ll}
\zeta_{k}&:=\disp\sum^{k}_{m=1}h_{k}\dist\big(w^{k}_{m},F(y^{k}_{m},t_{m})\big)=\disp\sum^{k}_{m=1}\disp\int^{t_m}_{t_{m-1}}\dist(w^{k}_{m}, F(y^{k}_{m},t_{m}))\\
&\le\disp\sum^{k}_{m=1}\disp\int^{t_m}_{t_{m-1}}\dist\big(w^{k}_{m},F(y^{k}_{m},t)\big)+\tau(F;h_{k}),\quad k\in\N.
\end{array}
\end{eqnarray}
It is well known (see, e.g., \cite[Proposition~6.3]{bs2}) that the a.e.\ continuity of $F(x,\cdot)$ on $[0,T]$ uniformly in $x\in U$ assumed in (H3) is equivalent to the convergence $\tau(F;h_{k})\to 0$.

Let us show next that $\disp\sum^{k}_{m=1}\disp\int^{t_m}_{t_{m-1}}\dist(w^{k}_{m},F(y^{k}_{m},t))\to 0$ as $k\to\infty$. Taking into account that $\bar{x}(\cdot)$ is a feasible trajectory for the differential inclusion (\ref{di}) and that $w^{k}(\cdot)\to\dot{\bar{x}}(\cdot)$ strongly in $L^1[0,T]$ and remembering also that for each $k\in\N$ the functions $w^{k}(t)$ are constant on the intervals $(t_{m-1},t_m]$, $m=1,2,\ldots$, and that $\dot{\ox}(t)$ is Riemann integrable (that is, a.e.\ continuous) on $[0,T]$, we can find $\tilde{t}_m\in(t_{m-1},t_m]$ such that
$$
\dot{\bar x}(\tilde{t}_m)\in F(\bar{x}(\tilde{t}_m),\tilde{t}_m)\;\mbox{ and }\;\disp\sum_{m=1}^k\int^T_0|\dot{\bar x}(\tilde{t}_m)-w^k(t)|\,dt\le 2\xi_k.
$$
This readily leads us to the following inequalities:
\begin{eqnarray*}
\begin{array}{lll}
&\disp\sum^{k}_{m=1}\disp\int^{t_m}_{t_{m-1}}\dist(w^{k}_{m},F(y^{k}_{m},t))dt\\
&\le\disp\sum^{k}_{m=1}\disp\int^{t_m}_{t_{m-1}}\Big[\dist(w^{k}_{m},F(\bar{x}_{m},t))+\dist(F(\bar{x}_{m},t),F(y^{k}_{m},t))\Big]dt\\
&\le\disp\sum^{k}_{m=1}\disp\int^{t_m}_{t_{m-1}}\Big[\dist(w^k_m,F(\bar{x}(\tilde{t}_m),\tilde{t}_m))+\dist(F(\bar{x}_{m},t),F(\bar{x}(\tilde{t}_m),\tilde{t}_m))\\
&\qquad\qquad\qquad+\dist(F(\bar{x}_m,t),F(y^{k}_{m},t))\Big]dt\\
&\le\disp\sum^{k}_{m=1}\int^{t_m}_{t_{m-1}}\Big[|w^k(t)-\dot{\bar x}(\tilde{t}_m)|+\dist(F(\bar{x}_{m},t),F(\bar{x}(\tilde{t}_m),\tilde{t}_m))+\dist(F(\bar{x}_{m},t),F(y^{k}_{m},t))\Big]dt\\
&\le2\xi_k+\disp\sum^{k}_{m=1}\int^{t_m}_{t_{m-1}}\Big[\dist(F(\bar{x}_{m},t),F(\bar{x}(\tilde{t}_m),\tilde{t}_m))+\dist(F(\bar{x}_{m},t),F(y^{k}_{m},t))\Big]dt\\
&\le2\xi_k+\tau(F;h_{k})+\disp\sum^{k}_{m=1}\int^{t_m}_{t_{m-1}}\Big[\dist(F(\bar{x}(t_m),t_m),F(\bar{x}(\tilde{t}_m),\tilde{t}_m))+\dist(F(\bar{x}_{m},t),F(y^{k}_{m},t))\Big]dt.
\end{array}
\end{eqnarray*}
Under the assumption (H3) we have $\disp\sum^{k}_{m=1}\disp\int^{t_m}_{t_{m-1}}\dist(w^{k}_{m},F(y^{k}_{m},t)) \to 0$. By employing (\ref{xi}) and the definition of $\zeta_k$ in (\ref{zeta}), it gives us the convergence $\zeta_k\to 0$ as $k\to\infty$.

Using this and the last inequality in (\ref{zy}) allows us to conclude that
\begin{eqnarray}\label{zjyjbound}
|z^{k}_{j+1}-y^{k}_{j+1}|\le\zeta_{k}e^{2lT}\;\mbox{ for all }\;j=0,\ldots,k-1\;\mbox{ and all }\;k\in\N.
\end{eqnarray}
Furthermore, we easily get the estimates
\begin{eqnarray}\label{eta}
|z^{k}_{j+1}-\bar{x}_{j+1}|\le\zeta_{k}e^{2lT}+|y^{k}_{j+1}-\bar{x}_{j+1}|\le\zeta_{k}e^{2lT}+\xi_k=:\eta_k,
\end{eqnarray}
where $\eta_k\to 0$ due to (\ref{xi}) and $\zeta_k\to 0$ as $k\to\infty$.

Considering next the the piecewise linear functions $z^k(\cdot)$ built in (\ref{dis-tr}) by using the discrete velocity $v^k(\cdot)$ from (\ref{dis-vel}), we get from (\ref{dis-vel}) and (\ref{iez}) that
$$
\dot{z}^{k}(t)=v^{k}_{j}=\frac{z^{k}_{j}-z^{k}_{j-1}}{h_k}\in F\big(z^k(t_j),t_j\big)\;\mbox{ on }\;(t_{j-1},t_j],\quad j=1,\ldots,k.
$$
It follows from the uniform boundedness of $F$ in (H1) that there is a subsequence of $\{\dot{z}^k(\cdot)\}$ (without relabeling) that converges to some function in $L^1[0,T]$, which cannot be anything but $\dot{\ox}(t)$ due to the relationships in (\ref{eta}) established above. Thus $\dot{z}^k(\cdot)\to\dot{\ox}(\cdot)$ {\em weakly} in $L^1[0,T]$ as $k\to\infty$. The latter is equivalent to the {\em uniform convergence} $z^k(\cdot)\to\ox(t)$ by the Newton-Leibniz formula (\ref{dis-tr}), and so we get (i).\vspace*{0.05in}

Now we justify assertion (ii) proving that in fact $\dot{z}^{k}(\cdot)\to\dot{\bar{x}}(\cdot)$ as $k\to\infty$ {\em strongly} in $L^{1}[0,T]$ provided that either the graph of $F(\cdot,t)$ is convex around $(\ox(t),\dot{\ox}(t))$, or the mapping $x\mapsto F(x,t)$ is locally Lipschitzian around $\ox(t)$ for a.e.\ $t\in[0,T]$.

First we examine the case when the {\em graph} of $F(\cdot,t)$ is {\em locally convex}. The classical Mazur's weak closure theorem tells us that there is a sequence of {\em convex combinations} of $\dot{z}^k(\cdot)$, which converges to $\dot{\ox}(\cdot)$ in the norm topology of $L^1[0,T]$ and thus contains a subsequence (no relabeling) converging to $\dot{\ox}(t)$ for a.e.\ $t\in[0,T]$. Taking into account the graph convexity of $F(\cdot,t)$ and the piecewise constant nature of $\dot{z}^k(t)$, we conclude that all the elements of the aforementioned sequence of convex combinations are {\em feasible} trajectories of the discrete approximation systems for any $k\in\N$. Therefore we get a sequence of feasible solutions to the discrete inclusions (\ref{ie}) whose piecewise linear extensions on $[0,T]$ converges to $\dot{\ox}(\cdot)$ strongly in $L^1[0,T]$. Keeping for simplicity the notation $\dot{z}^k(\cdot)$ for the elements of this sequence allows us to write
\begin{eqnarray}\label{al}
\al_k:=\disp\int^T_0|\dot{z}^k(t)-\dot{\ox}(t)|\,dt\to 0\;\mbox{ as }\;k\to\infty.
\end{eqnarray}

To complete the proof of the theorem in the convex graph case, it remains to verify the the convergence of $\{\dot{z}^{k}(\cdot)\}$ to $\dot{\ox}(\cdot)$ in the norm topology of $L^2[0,T]$. By the constructions above and assumption (H1), it is implied by the following relationships:
\begin{eqnarray}\label{al1}
\begin{array}{lll}
\disp\sum^{k}_{j=1}\disp\int^{t_{j}}_{t_{j-1}}\Big|\frac{z^{k}(t_{j})-z^{k}(t_{j-1})}{h_k}-\dot{\bar x}(t)\Big|^{2}dt&=\disp\sum^{k}_{j=1}\max\Big(|v^{k}_{j}|+|\dot{\bar x}(t)|\Big)\int^{t_{j}}_{t_{j-1}}|v^{k}_{j}-\dot{\bar x}(t)|\,dt\\
&\le 2m_{F}\disp\sum^{k}_{j=1}\disp\int^{t_{j}}_{t_{j-1}}|v^{k}_{j}-\dot{\bar{x}}(t)|\,dt=2m_{F}\alpha_k,
\end{array}
\end{eqnarray}
where $\al_k$ is taken from (\ref{al}). This justifies the $W^{1,2}[0,T]$-norm convergence of the extended discrete trajectories $\{z^k(\cdot)\}$ from (\ref{dis-tr}) in the first case under consideration.

Let us finally consider the other case in (ii) when $F(\cdot,t)$ is {\em locally Lipschitzian} around $\ox(t)$ for a.e.\ $t\in[0,T]$. Then for all $j\in\{1,\ldots,k-1\}$ we have the estimates
\begin{eqnarray}\label{ziyil}
\begin{array}{ll}
|z^{k}_{j+1}-z^{k}_{j}-h_{k}w_{j+1}|&\le\disp\frac{1}{1-lh_k}\dist(z^{k}_{j}+h_{k}w_{j+1},z^{k}_{j}+h_{k}F(z^{k}_{j}+h_{k}w_{j+1},t_{j+1}))\\
&\le\disp\frac{1}{1-lh_k}\dist(h_{k}w_{j+1},h_{k}F(z^{k}_{j}+h_{k}w_{j+1},t_{j+1}))\\
&\le\disp\frac{h_{k}}{1-lh_k}\Big[\dist(w_{j+1},F(y^{k}_{j}+h_{k}w_{j+1},t_{j+1}))\\
&+\dist(F(y^{k}_{j}+h_{k}w_{j+1},t_{j+1}),F(z^{k}_{j}+h_{k}w_{j+1},t_{j+1}))\Big]\\
&\le\disp\frac{h_{k}}{1-lh_k}\Big[l|z^{k}_{j}-y^{k}_{j}|+\dist(w_{j+1},F(y_{j+1},t_{j+1}))\Big].
\end{array}
\end{eqnarray}
Combining (\ref{ziyil}) with (\ref{zjyjbound}) gives us
\begin{eqnarray*}
\begin{array}{ll}
\disp\int^{T}_{0}|\dot{z}^k(t)-\dot{y}^k(t)|\,dt&=\disp\sum^{k}_{j=1}h_{k}|v^{k}_{j}-w^{k}_{j}|\\
&\le\disp\frac{1}{1-l h_{k}}\Big[lh_{k}\sum^{k}_{j=1}|z^{k}_{j-1}-y^{k}_{j-1}|+\disp\sum^{k}_{j=1}h_{k}\dist(w^{k}_{j},F(y^{k}_{j},t_{j}))\Big]\\
&\le 2(lT\zeta_{k}e^{2lT}+\zeta_{k})=2\zeta_{k}(lTe^{2lT}+1)\to 0\;\mbox{ as }\;k\to\infty,
\end{array}
\end{eqnarray*}
where $\zeta_k\to 0$ is taken from (\ref{zeta}). Taking further into account that
\begin{eqnarray*}
\disp\int^{T}_{0}|\dot{z}^{k}(t)-\dot{\bar x}(t)|\,dt\le\int^{T}_{0}|\dot{z}^{k}(t)-\dot{y}^k(t)|\,dt+\int^{T}_{0}|\dot{y}^{k}(t)-\dot{\bar x}(t)|\,dt
\end{eqnarray*}
and using the convergence $\xi_k\to 0$ in (\ref{xi}) with $\dot{y}^k(t)=w^k(t)$ tell us that the $L^1[0,1]$-norm convergence of $\{\dot{z}^k(\cdot)\}$  in (\ref{al}) holds in the second case under consideration. Applying now (\ref{al1}) justifies (ii) in this case and thus completes proof of the theorem. $\h$

\section{Strong Convergence of Discrete Optimal Solutions}\sce

In this section we construct a sequence of well-posed discrete approximations of the Bolza problem $(P)$ for ROSL differential inclusions and
justify the norm convergence in either $C[0,T]$ or $W^{1,2}[0,T]$ topology of their optimal solutions to either a strong local minimizer or an intermediate relaxed local minimizer $\ox(\cdot)$ of $(P)$, respectively. In addition to our standing assumptions (H1)--(H3) on the right-hand side $F$ in (\ref{di}) and those (if needed) from Theorem~\ref{sa} formulated now around the given local minimizer, the following ones are imposed here on the functions $f$ and $\ph_0$ in the Bolza cost functional (\ref{cost}) as well as on the functions $\ph_i$, $i=1,\ldots,m+r$, and the set $\O$ in the endpoint constraints (\ref{edc})--(\ref{equ}).

\begin{enumerate}
\item[(\bf H4)]\label{h3}
The function $f(x,v,\cdot)$ is a.e.\ continuous on $[0,T]$ and bounded uniformly in $(x,v)\in U\times(m_{F}B)$. Furthermore,
there exists $\nu>0$ such that the function $f(\cdot,\cdot,t)$ is continuous on the set
\begin{eqnarray*}
A_{\nu}(t)=\big\{(x,v)\in U\times(m_F+\nu)\B\big|\;v\in F(x,t')\;\mbox{ for some }\;t'\in(t-\nu,t]\big\}
\end{eqnarray*}
uniformly in $t$ on the interval $[0,T]$.

\item[(\bf H5)]\label{h4}
The cost function $\varphi_0$ is continuous on $U$, while the constraint functions $\varphi_i$ are Lipschitz continuous on $U$ for all $i=1,\ldots,m+r$. Furthermore, the endpoint constraint set $\Omega$ is locally closed around $\ox(T)$.
\end{enumerate}

Given a {\em r.i.l.m.}\ $\bar{x}(\cdot)$ in ($P$), suppose without loss of generality (due to (H1)) that $\al=1$ and $p=2$ in (\ref{ilme}) and Definition~\ref{rilm}. Denote by $L>0$ a common Lipschitz constant for the functions $\ph_i$, $i=1,\ldots,m+r$, on $U$ and take the sequence $\{\eta_k\}$ in (\ref{eta}) constructed via the approximation of the local optimal solution $\ox(\cdot)$ under consideration. Then we define a sequence of discrete approximation problems $(P_k)$, $k\in\N$, as follows:
\begin{eqnarray}\label{minpk}
\begin{array}{lll}
\mbox{minimize}\;J_k[x^k]:&=\varphi_0(x^k(t_k))+h_{k}\disp\sum^{k}_{j=1}f\Big(x^{k}(t_{j}),\frac{x^{k}(t_{j})-x^{k}(t_{j-1})}{h_k},t_{j}\Big)\\
&+\disp\sum^{k}_{j=1}\disp\int^{t_{j}}_{t_{j-1}}\Big|\frac{x^{k}(t_{j})-x^{k}(t_{j-1})}{h_k}-\dot{\bar x}(t)\Big|^{2}dt
\end{array}
\end{eqnarray}
over trajectories $x^k=(x^k_0,\ldots,x^k_k)$ of the discrete inclusions (\ref{ie}) subject to the constraints
\begin{equation}\label{pkxk}
|x^{k}(t_j)-\bar{x}(t_j)|^2\le\frac{\varepsilon^2}{4}~\;\mbox{for}~\;j=1,\ldots,k,
\end{equation}
\begin{equation}\label{pkd}
\sum^{k}_{j=1}\int^{t_{j}}_{t_{j-1}}\Big|\frac{x^{k}(t_{j})-x^{k}(t_{j-1})}{h_k}-\dot{\bar x}(t)\Big|^{2}dt\le\frac{\varepsilon}{2},
\end{equation}
\begin{equation}\label{edcg}
x^k_k\in\Omega_{k}:=\Omega+\eta_{k}\B,
\end{equation}
\begin{equation}\label{le0}
\varphi_{i}(x^k_k)\le L\eta_k\;\mbox{ for }\;i=1,\ldots,m,
\end{equation}
\begin{equation}\label{eq0}
-L\eta_k\le\varphi_{i}(x^k_k)\le L\eta_k\;\mbox{ for }\;i=m+1,\ldots,m+r,
\end{equation}
where $\ve>0$ is fixed and taken from (\ref{ilme}) for the given r.i.l.m.\ $\ox(\cdot)$.

If $\ox(\cdot)$ is a given {\em strong local minimizer} for $(P)$ with $f=f(x,t)$, we construct a simplified sequence of discrete approximations problems $(\tilde P_k)$ as follows:
\begin{eqnarray}\label{minpk1}
\mbox{minimize}\;\tilde J_k[x^k]:&=\varphi_0(x^k(t_k))+h_{k}\disp\sum^{k}_{j=1}f(x^{k}(t_{j}),t_{j})+\disp\sum^{k}_{j=1}|x^k(t_j)-\ox(t_j)|^{2}
\end{eqnarray}
subject to the constraints (\ref{pkxk})--(\ref{eq0}) with $\eta_k$ taken from (\ref{eta}).

The next theorem shows that problems $(P_k)$ and $(\tilde P_k)$ admit optimal solutions for all large $k\in\N$ and that extended discrete optimal solutions to these problems converge to $\ox(\cdot)$ in the corresponding norm topology of either $C[0,T]$ or $W^{1,2}[0,T]$ depending on the type of local minima (strong or intermediate) which $(P)$ achieves at $\ox(\cdot)$.

\begin{Theorem}[strong convergence of discrete optimal solutions]\label{dsa} Let $\bar x(\cdot)$ be a Riemann integrable local optimal solution to the original Bolza problem $(P)$ under the validity of assumptions {\em (H1)-(H5)} around $\ox(\cdot)$. The following assertions hold:

{\bf (i)} If $\ox(\cdot)$ is a strong local minimizer for $(P)$ with $f=f(x,t)$, then each problem $(\Tilde P_k)$ admits an optimal solution $\ox^k(\cdot)$ for large $k\in\N$ and the sequence $\{\bar{x}^{k}(\cdot)\}$ piecewise linearly extended to $[0,T]$ converges to $\bar{x}(\cdot)$ as $k\to\infty$ in the norm topology of $C[0,T]$.

{\bf (ii)} If $\ox(\cdot)$ is a r.i.l.m.\ in $(P)$ and the assumptions of Theorem~{\rm\ref{sa}(ii)} are satisfied for $\ox(\cdot)$, then each problem $(P_k)$ admits an optimal solution $\ox^k(\cdot)$ whenever $k\in\N$ is sufficiently large and the sequence $\{\bar{x}^{k}(\cdot)\}$ piecewise linearly extended to $[0,T]$ converges to $\bar{x}(\cdot)$ as $k\to\infty$ in the norm topology of $W^{1,2}[0,T]$.
\end{Theorem}
{\bf Proof.} We verify the existence of optimal solutions to problems $(P_k)$ and $(\tilde P_k)$ in a parallel way. Observe first that both $(P_k)$ and $(\tilde P_k)$ admit feasible solutions for all $k\in\N$ sufficiently large. Indeed, take for each $k$ the discrete trajectories $z^k:=(z^k_0,\ldots,z^k_k)$ constructed in Theorem~\ref{sa}(i) to approximate the r.i.l.m.\ $\ox(\cdot)$ and in Theorem~\ref{sa}(ii) to approximate the r.i.l.m.\ $\ox(\cdot)$. Then both these functions satisfy the discrete inclusion (\ref{ie}), and it remains to verify that the corresponding $z^k$  fulfills the constraints in (\ref{pkxk}), (\ref{edcg})--(\ref{eq0}) in the case of $(\tilde P_k)$ and those in (\ref{pkxk})--(\ref{eq0}) in the case of $(P_k)$. The validity of (\ref{pkxk}) and (\ref{edcg}) in both cases follows from (\ref{eta}) for large $k$, while the validity of the additional constraint (\ref{pkd}) for $(P_k)$ follows from (\ref{al}). The fulfillment of the inequality constraints in (\ref{le0}) and (\ref{eq0}) for $z^k_k$ follows by these arguments from the validity of (\ref{ine}) and (\ref{equ}) for $\ox(T)$, respectively, and the local Lipschitz continuity of the endpoint functions
$$
|\varphi_{i}(z^{k}_{k})-\varphi_{i}(\ox(T))|\le L|z^{k}_{k}-\ox(T)|\le L\eta_{k},\quad i=1,\ldots,m+r.
$$
Thus for each $k\in\N$ (omitting the expression ``for all large $k$" in what follows) the sets of feasible solutions to $(P_k)$ and $(\tilde P_k)$ are nonempty. It is clear from the construction of $(P_k)$ and $(\tilde P_k)$ and the assumptions made that each of these sets is closed and bounded. This ensures the existence of optimal solutions to $(P_k)$ by the classical Weierstrass existence theorem due to the continuity of the functions $\ph_0$ and $f$ in (\ref{minpk}) and (\ref{minpk1}) .\vspace*{0.05in}

Next we proceed with the proof of the strong $W^{1,2}[0,T]$-convergence in (ii) for any sequence of the discrete optimal solutions $\{\ox^k(\cdot))\}$ in $(P_k)$ piecewise linearly extended to the continuous-time interval $[0,T]$. To this end let us first show that
\begin{equation}\label{liminf}
\disp\liminf_{k\to\infty}J_{k}[\bar{x}^k]\le J[\bar x]
\end{equation}
for the optimal values of the cost functional in (\ref{minpk}). It follows from the optimality of $\ox^k(\cdot)$ for $(P_k)$ and the feasibility of $z^k(\cdot)$ taken from the proof of (i) for this problem that $J_{k}[\bar{x}^k]\le J_{k}[z^k]$ for each $k$. To get (\ref{liminf}), it suffices to show therefore that
\begin{equation}\label{lim1}
J_{k}[z^k]\to J[\bar x]\;\mbox{ as }\;k\to\infty
\end{equation}
including the verification of the existence of the limit. We have from (\ref{minpk}) that
\begin{eqnarray*}
J_k[z^k]&=\varphi_0(z^k(t_k))
+h_{k}\disp\sum^{k}_{j=1}f\Big(z^{k}(t_{j}),\frac{z^{k}(t_{j})-z^{k}(t_{j-1})}{h_k},t_{j}\Big)\\
&+\disp\sum^{k}_{j=1}\disp\int^{t_{j}}_{t_{j-1}}\Big|\frac{z^{k}(t_{j})-z^{k}(t_{j-1})}{h_k}-\dot{\bar x}(t)\Big|^{2}dt
\end{eqnarray*}
and deduce from the convergence $z^k(t_k)\to\ox(T)$ and the continuity assumption on $\ph_0$ in (H5) the convergence $\varphi_0(z^k(t_k))\to\varphi_0(\bar{x}(T))$ as $k\to\infty$ of the terminal cost function in (\ref{minpk}). Furthermore, it follows from (\ref{al}) that
$$
\disp\sum^{k}_{j=1}\disp\int^{t_{j}}_{t_{j-1}}\Big|\frac{z^{k}(t_{j})-z^{k}(t_{j-1})}{h_k}-\dot{\bar x}(t)\Big|^{2}dt\to 0\;\mbox{as }\;k\to\infty.
$$
To justify (\ref{lim1}), we only need to check that
$$
h_{k}\sum^{k}_{j=1}f\Big(z^{k}(t_{j}),\frac{z^{k}(t_{j})-z^{k}(t_{j-1})}{h_k},t_{j}\Big)\to\int^{T}_{0}f(\bar{x}(t),\dot{\bar{x}}(t),t)dt~\;\mbox{as}~\;k \to \infty.
$$
The continuity assumptions on $f$ in (H4) imply without loss of generality that
$$\Big|f(z^{k}(t_{j}),\frac{z^{k}(t_{j})-z^{k}(t_{j-1})}{h_k},t_{j})-f(z^{k}(t_{j}),\frac{z^{k}(t_{j})-z_{k}(t_{j-1})}{h_k},t)\Big|\le\frac{\varepsilon}{T}
$$
for all $k\in\N$ and a.e.\ $t\in[0,T]$. Employing now Lebesgue's dominated convergence theorem together with Theorem~\ref{sa}(ii) tells us that
\begin{eqnarray*}
\begin{array}{ll}
& h_{k}\disp\sum^{k}_{j=1}f\Big(z^{k}(t_{j}),\frac{z^{k}(t_{j})-z^{k}(t_{j-1})}{h_k},t_{j}\Big)=\disp\sum^{k}_{j=1}\disp\int^{t_{j}}_{t_{j-1}}f(z^{k}(t_{j}),v^{k}(t),
t_{j})\,dt\\&\sim\disp\sum^{k}_{j=1}\disp\int^{t_{j}}_{t_{j-1}}f(z^{k}(t_{j}),v^{k}(t),t)\,dt\sim\disp\sum^{k}_{j=1}\disp\int^{t_{j}}_{t_{j-1}}f(\bar{x}(t),v^{k}(t),t)\,dt\\
&=\disp\int^{T}_{0}f(\bar{x}(t),v^{k}(t),t)\,dt\sim\disp\int^{T}_{0}f(\bar{x}(t),\dot{\bar{x}}(t),t)\,dt,
\end{array}
\end{eqnarray*}
where the sign '$\sim$' is used to indicate the equivalence as $k\to\infty$. Thus we get (\ref{lim1}).

To proceed further, consider the numerical sequence
\begin{eqnarray}\label{c}
c_{k}:=\disp\int^{T}_{0}|\dot{\bar{x}}^k(t)-\dot{\bar{x}}(t)|^2dt,\quad k\in\N,
\end{eqnarray}
and verify that $c_{k}\to 0$ as $k\to\infty$. Since the numerical sequence in (\ref{c}) is obviously bounded, it has limiting points. Denote by $c\ge 0$ any of them and show that $c=0$. Arguing by contradiction, suppose that $c>0$. It follows from the uniform boundedness assumption (H1) and basic functional analysis that the sequence $\{\dot{\ox}^k(\cdot)\}$ contains a subsequence (without relabeling), which converges in the weak topology of $L^2[0,T]$ to some $v(\cdot)\in L^2[0,T]$. Considering the absolutely continuous function
$$
\widetilde{x}(t):=x_0+\disp\int^{t}_{0}v(s)\,ds,\quad 0\le t\le T,
$$
we deduce from the Newton-Leibniz formula that the sequence of the extended discrete trajectories $\ox^k(\cdot)$ converges to $\Tilde x(\cdot)$ in the weak topology of $W^{1,2}[0,T]$, for which we have $\dot{\Tilde x}(t)=v(t)$ for a.e.\ $t\in[0,T]$. By invoking Mazur's weak closure theorem, it follows from the convexity of the sets $F(x,t)$ and the continuity of $F(\cdot,t)$ that the limiting function $\Tilde x(\cdot)$ satisfies the differential inclusion (\ref{di}). Furthermore, the construction of the discrete approximation problems $(P_k)$ with $\eta_k\to 0$ therein ensures that $\Tilde x(\cdot)$ is a feasible trajectory for the original Bolza problem $(P)$, and therefore for the relaxed $(R)$ as well.

Employing again Mazur's weak closure theorem, we find a sequence of convex combination of $\dot{\bar x}^k(\cdot)$ converging to $\dot{\widetilde x}(\cdot)$ in the norm topology of $L^2[0,T]$ and hence a.e.\ on $[0,T]$ along some subsequence. Taking into account the construction of $\widehat{f}_{F}$ as the convexification of $f_F$ in (\ref{ext-f}) with respect to the velocity variable, we arrive at the inequality
\begin{equation}\label{relax-p}
\int^{T}_{0}\widehat{f}_{F}(\widetilde{x}(t),\dot{\widetilde x}(t),t)\,dt\le\disp\liminf_{k\to\infty}h_{k} \sum^{k}_{j=1}f\Big(\bar{x}^{k}_{j},\frac{\bar{x}^{k}_{j}-\bar{x}^{k}_{j-1}}{h_{k}},t_{j}\Big).
\end{equation}

Define now the integral functional on $L^2[0,T]$ by
\begin{eqnarray}\label{I}
I[v]:=\disp\int^{T}_{0}|v(t)-\dot{\bar x}|^2dt
\end{eqnarray}
and show it is convex on this space. Indeed, picking any $v(\cdot),w(\cdot)\in L^2[0,T]$ and $\lm\in[0,1]$ and using the Cauchy-Schwartz inequality gives us
\begin{eqnarray*}
I[\lambda v+(1-\lambda)w]&=&\int^{T}_{0}|\lambda(v(t)-\dot{\bar x}(t))+(1-\lambda)(w(t)-\dot{\bar x}(t))|^2dt\\
&\le&\int^{T}_{0}\Big[\lm|v(t)-\dot{\bar x}(t)|+(1-\lambda)|w(t)-\dot{\bar x}(t)|\Big]^2dt\\
&=&\lambda\int^{T}_{0}|v(t)-\dot{\bar x}(t)|^2dt+(1-\lambda)\int^{T}_{0}|w(t)-\dot{\bar x}(t)|^2dt\\
&=&\lambda I[v]+(1-\lambda)I[w],
\end{eqnarray*}
which justifies the convexity and hence the lower semicontinuity of (\ref{I}) in the weak topology of $L^2[0,T]$. It allows us to conclude that
\begin{eqnarray*}
\int^{T}_{0}|\dot{\widetilde x}(t)-\dot{\bar x}(t)|^{2}dt&\le&\disp\liminf_{k\to\infty}\int^{T}_{0}|\dot{\bar x}^k(t)-\dot{\bar x}(t)|^{2}dt\\
&=&\disp\liminf_{k\to\infty}\sum^{k}_{j=1}\int^{t_{j}}_{t_{j-1}}\Big|\frac{\bar{x}^{k}(t_{j})-\bar{x}^{k}(t_{j-1})}{h_k}-\dot{\bar x}(t)\Big|^2dt.
\end{eqnarray*}
Employing this and passing to the limit in the constraints (\ref{pkxk}) and (\ref{pkd}) for $\ox^k(\cdot)$ yield
\begin{eqnarray*}
|\widetilde{x}(t)-\bar{x}(t)|\le\disp\frac{\varepsilon}{2}\;\mbox{ for }\;t\in[0,T]\;\mbox{ and }\;\disp\int^{T}_{0}|\dot{\widetilde x}(t)-\dot{\bar x}(t)|^{2}dt\le\disp\frac{\varepsilon}{2},
\end{eqnarray*}
which verifies that the feasible trajectory $\Tilde x(\cdot)$ for $(R)$ belongs to the prescribed $W^{1,2}[0,T]$ neighborhood of the r.i.l.m.\ $\ox(\cdot)$ from Definition~\ref{rilm}.

Now we are able to pass to the limit in the cost functional formula (\ref{minpk}) in $(P_k)$ for $\ox^k(\cdot)$ by using (\ref{liminf}), (\ref{relax-p}), and the assumption on $c_k\to c>0$ in (\ref{c}). It gives us
$$
\Hat J[\Tilde x]=\varphi(\widetilde{x}(T))+\disp\int^{T}_{0}\widehat{f}_{F}(\widetilde{x}(t),\dot{\widetilde x}(t),t)\,dt\le\liminf_{k\to\infty}J_{k}[\bar{x}^{k}]+c<J[\ox]=\Hat J[\ox],
$$
which contradicts the choice of $\bar x(\cdot)$ as a r.i.l.m.\ for the original Bolza problem $(P)$. Thus we have $c_k\to 0 $ as $k\to\infty$ showing in this way that $\ox^k(\cdot)\to\ox(\cdot)$ strongly in $W^{1,2}[0,T]$.\vspace*{0.05in}

To complete the proof of the theorem, it remains to justify the strong $C[0,T]$ convergence in (i) of discrete optimal trajectories for $(\tilde P_k)$ in the case when $\ox(\cdot)$ is a strong local minimizers of the continuous-time Bolza problem $(P)$. Note that due to the convexity of $F(x,t)$ and the independence of the integrand $f$ on the velocity variable, problem $(P)$ agrees with its relaxation $(R)$. Taking into account the form of the cost functional (\ref{minpk1}) and Theorem~\ref{sa}(i) on the strong discrete approximation of $\ox(\cdot)$ in $C[0,T]$, we arrive at the claimed convergence result in assertion (i) of this theorem by just simplifying the above proof of assertion (ii) and replacing the cost functional $J_k$ with $\tilde J_k$. $\h$

\section{Optimality Conditions for Discrete Approximations}\sce

In this section we derive necessary optimality conditions for each problem $(P_k)$, $k\in\N$, in the sequence of discrete approximations formulated in Section~4. In the same way we can proceed with deriving necessary optimality conditions in the simplified problems $(\tilde P_k)$; we do not present them here due to the full similarity and size limitation.

Note that problems of this type intrinsically belong to nonsmooth optimization even when all the functions $f$ and $\ph_i$ for $i=0,\ldots,m+r$ are smooth and $\O=\R^n$. The nonsmoothness comes from the dynamic constraints in (\ref{ie}) given by the discretization of the differential inclusion (\ref{di}); the number of these constraints is increasing along with decreasing the step of discretization. To derive necessary optimality conditions for problems $(P_k)$, we employ advanced tools of variational analysis and generalized differentiation discussed in Section~2.\vspace*{0.05in}

Consider first the following problem of mathematical programming $(MP)$ with finitely many functional and geometric constraints. Given $\phi_{j}\colon\R^{d}\to\R$ for $j=0,\ldots,s$, $g_{j}:\R^{d}\to\R^{n}$ for $j=0,\ldots,p$, and $\Delta_j\subset\R^d$ for $j=0,\ldots,q$, we define ($MP$) by
\begin{eqnarray*}
\begin{array}{ll}
&\mbox{minimize }\;\phi_{0}(z)\;\mbox{ subject to}\\
&\phi_{j}(z)\le 0~\;\mbox{for}~\;j=0,\ldots,s,\\
&g_{j}(z)=0~\;\mbox{for}~\;j=0,\ldots,p,\\
&z\in\Delta_{j}~\;\mbox{for}~\;j=0,\ldots,q.
\end{array}
\end{eqnarray*}

The next result gives us necessary optimality conditions for local minimizers of problem $(MP)$ in the setting needed for the subsequent application to deriving optimality conditions in the discrete approximation problems $(P_k)$. We express these conditions via our basic normal cone (\ref{mnc}) and subdifferential (\ref{subd}) constructions from Section~2.

\begin{Lemma}{\bf(generalized Lagrange multiplier rule for mathematical programs).}\label{lmmp}
Let $\bar z$ be a local optimal solution to problem $(MP)$. Assume that the functions $\phi_{j}$ are Lipschitz continuous around $\oz$, the mappings $g_{j}$ are continuous differentiable around $\oz$, and the sets $\Delta_j$ are locally closed around this point. Then there exist nonnegative numbers $\mu_{j}$ for $j=0,\ldots,s$ as well as vectors $\psi_{j}\in\R^{n}$ for $j=0,\ldots,p$ and $z^{*}_{j}\in\R^{d}$ for $j=0,\ldots,q$, not equal to zero simultaneously, such that we have the conditions
\begin{eqnarray*}
z^{*}_{j}\in N(\bar z;\Delta_{j}),\quad j=0,\ldots,q,
\end{eqnarray*}
\begin{eqnarray*}\label{slack}
\mu_{j}\phi_{j}(\bar z)=0,\quad j=1,\ldots,s,
\end{eqnarray*}
\begin{eqnarray*}\label{lag}
-z^{*}_{0}-\ldots-z^{*}_{q}\in\partial\Big(\sum^{s}_{j=0}\mu_{j}\phi_{j}\Big)(\bar z)+\sum^{p}_{j=0}(\nabla g_{j}(\bar z))^{T}\psi _{j},
\end{eqnarray*}
where the symbol ``$A^T$" indicates the matrix transposition.
\end{Lemma}
{\bf Proof.} This result follows from necessary optimality conditions given \cite[Theorem~5.21]{bs2} for problems with a single geometric constraint and the basic intersection rule for the normal cone (\ref{mnc}) taken from \cite[Theorem~3.4]{bs1}. $\h$\vspace*{0.05in}

Now we employ Lemma~\ref{lmmp} and calculus rules for generalized normals and subgradients to derive necessary optimality conditions for the structural dynamic problems of discrete approximation $(P_k)$ in the extended Euler-Lagrange form. Note that for this purpose we need less assumptions that those imposed in (H1)--(H5). Observe also that the form of the Euler-Lagrange inclusion below reflects the essence of the implicit Euler scheme being significantly different from the adjoint system corresponding to the explicit Euler counterpart from \cite{bs95,bs2}. The solvability of the new implicit adjoint system is ensures by Lemma~\ref{lmmp} due the given proof of the this theorem.

\begin{Theorem}[extended Euler-Lagrange conditions for discrete approximations]\label{dnc}
Fix any $k\in\N$ and let $\bar x^{k}=(\bar x^{k}_{0},\ldots, \bar{x}^{k}_{k})$ with $\bar{x}^k_0=x_0$ in {\rm(\ref{ie})} be an optimal solution to problem $(P_{k})$ constructed in Section~{\rm 4}. Assume that the sets $\Omega$ and $\gph F_{j}$ with $F_j:=F(\cdot,t_j)$ are closed and the functions $\varphi_{i}$ for $i=0,\ldots,m+r$ and $f_{j}:=f(\cdot,\cdot,t_j)$ for $j=0,\ldots,k$ are Lipschitz continuous around the corresponding points.

Then there exist real numbers $\lambda^{k}_{i}$ for $i=0,\ldots,m+r$ and a vector $p^{k}:=(p^{k}_{0},\ldots, p^{k}_{k})\in\R^{(k+1)n}$, which are not equal to zero simultaneously and satisfy the following relationships:

$\bullet$ The sign conditions
$$
\lm^{k}_i\ge 0\;\mbox{ for }\;i=0,\ldots,m;
$$

$\bullet$ the complementary slackness conditions
$$
\lambda^{k}_{i}[\varphi_i(\bar{x}^k_k)-L\eta_k]=0\;\mbox{ for }\;i=1,\ldots,m;
$$

$\bullet$ the extended Euler-Lagrange inclusion held for $j=1,\ldots,k$:
\begin{eqnarray*}
\disp\Big(\frac{p^{k}_{j}-p^k_{j-1}}{h_k},p^{k}_{j-1}-\frac{\lambda^{k}_{0}\theta^{k}_{j}}{h_k}\Big)\in\lambda^{k}_{0}
\disp\partial f_j\Big(\bar{x}^{k}_{j},\frac{\bar{x}^{k}_{j}-\bar{x}^{k}_{j-1}}{h_k}\Big)+N\Big(\Big(\bar{x}^{k}_{j},
\frac{\ox^k_{j}-\ox^k_{j-1}}{h_k}\Big);\gph F_j\Big);
\end{eqnarray*}

$\bullet$ the transversality inclusion
\begin{eqnarray*}
-p^{k}_{k}\in\sum^{m}_{i=0}\lambda^{k}_{i}\partial\varphi_{i}(\bar{x}^{k}_{k})+\sum^{m+r}_{i=m+1}\lambda^{k}_{i}\partial^{0}\varphi_{i}(\bar{x}^{k}_{k})
+N(\bar{x}^{k}_{k};\Omega_k),
\end{eqnarray*}
where $\partial^0\ph_i$ stands for the symmetric subdifferential {\rm(\ref{sym})} of $\ph_i$, and where
\begin{eqnarray}\label{theta}
\disp\theta^{k}_{j}:=-2\int^{t_{j}}_{t_j-1}\Big(\dot{\bar x}(t)-\frac{\bar{x}^{k}_{j}-\bar{x}^{k}_{j-1}}{h_k}\Big)\,dt.
\end{eqnarray}
\end{Theorem}
{\bf Proof.} Skipping for notational simplicity the upper index ``$k$" if no confusions arise, consider the new ``long" variable
$$
z:=(x_{0},\ldots,x_{k},y_{1},\ldots,y_{k})\in\mathbb{R}^{(2k+1)n}\;\mbox{ with the fixed initial vector }\;x_0
$$
and for each $k\in\N$ reformulate the discrete approximation problem $(P_k)$ as a mathematical program of the above type $(MP)$ with the following data:
\begin{equation}\label{mp1}
\mbox{min}\;\phi_{0}(z):=\varphi_0(x_k)+h_{k}\sum^{k}_{j=1} f(x_{j}, y_{j}, t_{j})+\sum ^{k}_{j=1}\int^{t_{j}}_{t_j-1}|y_{j}-\dot{\bar x}(t)|^{2}dt
\end{equation}
subject to the functional and geometric constraints
\begin{equation}\label{phij}
\phi_{j}(z):=|x_{j}-\bar{x}(t_{j})|^2-\disp\frac{\varepsilon^2}{4}\le 0\;\mbox{ for }\;j=1,\ldots,k,
\end{equation}
\begin{equation}\label{mp2}
\phi_{k+1}(z):=\sum^{k}_{j=1}\int^{t_{j}}_{t_j-1}|y_{j}-\dot{\bar x}(t)|^{2}dt-\frac{\varepsilon}{2}\le 0,
\end{equation}
\begin{eqnarray}\label{phimr}
\phi_{k+1+j}(z)=\varphi_{j}(x_k)-L\eta_k\le 0\;\mbox{ for }\;j=1,\ldots,m+r,
\end{eqnarray}
\begin{eqnarray}\label{phir}
\phi_{k+1+m+r+j}(z):=-\varphi_{m+j}(x_k)-L\eta_k\le 0\;\mbox{ for }\;j=1,\ldots,r,
\end{eqnarray}
\begin{equation}\label{mp3}
g_{j}(z):=x_{j}-x_{j-1}-h_{k}y_{j}=0\;\mbox{ for }\;j=1,\ldots,k,\;\;g_0(z)=x(0)-x_0\equiv 0,
\end{equation}
\begin{equation}\label{delta0}
z\in\Delta_{0}=\{(x_{0},\ldots,x_{k},y_{1},\ldots,y_{k})\in\mathbb{R}^{(2k+1)n}|\;x_{k}\in\Omega\},
\end{equation}
\begin{equation}\label{delta}
z\in\Delta_{j}=\{(x_{0},\ldots, x_{k},y_{1},\ldots, y_{k})\in \mathbb{R}^{(2k+1)n}|\;y_{j}\in F_{j}(x_{j})\},\quad j=1,\ldots,k.
\end{equation}

Let $\bar{x}^{k}=(x_0,{\ox}^{k}_{1},\ldots,{\ox}^{k}_{k})$ be a given local optimal solution to problem $(P_{k})$, and thus the corresponding extended variable
$\bar{z}:=(x_{0},\ldots,{\ox}_{k},(\bar{x}_{1}-\bar{x}_{0})/h_{k},\ldots,(\bar{x}_{k}-\bar{x}_{k-1})/h_{k})$, where the upper index ``$k$" is omitted, gives
a local minimum to the mathematical program $(MP)$ with the data defined in (\ref{mp1})--(\ref{delta}). Applying now to $\oz$ the generalized Lagrange multiplier rule from Lemma~\ref{lmmp}, we find normal collections
\begin{eqnarray}\label{nor1}
z^*_{j}=({x}^*_{0j},\ldots,{x}^*_{kj},y^*_{1j},\ldots,y^*_{kj})\in N(\bar{z};\Delta_{j})\;\mbox{ for }\;j=0,\ldots,k
\end{eqnarray}
and well as nonnegative multipliers $(\mu_{0},\ldots,\mu_{k+1+m+2r})$ and vectors $\psi_{j}\in\mathbb{R}^{n}$ for $j=0,\ldots,k$
such that we have the conditions
\begin{eqnarray}\label{slack}
\mu_{j}\phi_{j}(\bar z)=0\;\mbox{ for }\;j=1,\ldots,k+1+m+2r,
\end{eqnarray}
\begin{eqnarray}\label{lag}
-z^*_{0}-\ldots-z^*_{k}\in\partial\Big(\sum^{k+1+m+2r}_{j=0}\mu_{j}\phi_{j}\Big)(\bar z)+\sum^{k}_{j=0}(\nabla g_{j}(\bar z))^{T}\psi_{j}.
\end{eqnarray}
It follows from (\ref{nor1}) and the structure of $\Delta_0$ in (\ref{delta0}) that
$$
x^*_{k0}\in N(\bar{x}_{k};\Omega_{k}),\;y^*_{i0}=0\;\mbox{ for }i=0,\ldots,k-1,\;x^*_{i0}=0\;\mbox{ for }\;i=1,\ldots,k-1,\;\mbox{ and }\;x^*_{00}\;\mbox{ is free};
$$
the latter is due to the fact that $x_0$ is fixed. Furthermore, inclusion (\ref{nor1}) for $j=1,\ldots,k$ gives us by the structure of $\Delta_j$ that
$$
(x^*_{jj},y^*_{jj})\in N\Big(\Big(\bar{x}_{j},\disp\frac{\bar {x}_{j}-\bar {x}_{j-1}}{h_k}\Big);\gph F_{j}\Big)\;\mbox{ and }\;x^*_{ij}=y^*_{ij}=0\;\mbox {if }\;i\ne j,\;j=1,\ldots,k.
$$
Employing the above conditions together with the subdifferential sum rule from \cite[Theorem~2.33]{bs1} with taking into the nonnegativity of $\mu_j$, we get from (\ref{lag}) that
\begin{eqnarray*}
\begin{array}{ll}
&\partial\disp\Big(\sum^{k+1+m+2r}_{j=0}\mu_{j}\phi_{j}\Big)(\bar z)+\disp\sum^{k}_{j=0}(\nabla g_{j}(\bar z))^T\psi_{j}\subset \disp\sum^{k+1+m+2r}_{j=0}\mu_{j}\partial\phi_{j}(\bar z)+\disp\sum^{k}_{j=0}(\nabla g_{j}(\bar z))^T\psi_{j}\\
&\;\;=\mu_{0}\nabla\Big[\varphi(x_k)+h_{k}\disp\sum^{k}_{j=1}f(x_{j},y_{j}, t_{j})+\disp\sum^{k}_{j=1}\disp\int^{t_{j}}_{t_j-1}|y_{j}-\dot{\bar x}(t)|^{2}dt\Big]+\disp\sum^{k}_{j=1}\mu_{j}\nabla(|x_{j}-\bar{x}(t_{j})|^2)\\
&\;\;+\mu_{k+1}\nabla\Big(\disp\sum^{k}_{j=1}\int^{t_{j}}_{t_j-1}|y_{j}-\dot{\bar x}(t)|^{2}dt\Big)+\disp\sum^{m+r}_{j=1}\mu_{k+1+j}\nabla\varphi_{i}(\ox_k)-\disp\sum^{r}_{j=1}\mu_{k+1+m+r+j}\nabla\varphi_{j}(\ox_k)\\
&\;\;+\disp\sum^{k}_{j=1}\nabla(x_{j}-x_{j-1}-h_{k}y_{j})^T\psi_{j}+\nabla(x(0)-x_0)^T\psi_{0},
\end{array}
\end{eqnarray*}
where the derivatives (gradients, Jacobians) of all the {\em composite/sum} functions involves with respect of all their variables of are taken at the optimal point $\oz$.  It follows from Theorem~\ref{sa} that for $k\in\N$ sufficiently large we have $\phi_{j}(\bar{z}^{k})<0$ for $\oz=\oz^k$ and $j=1,\ldots,k+1$ due to the structures of the functions $\phi$ in (\ref{phij}) and (\ref{mp2}) and the complementary slackness conditions in (\ref{slack}). This implies $\mu_{j}=0$ for $j=1,\ldots,k+1$. Considering now the Lagrange multipliers
$$
\lambda^{k}_{0}:=\mu_{0}\;\mbox{ and }\;\lambda^{k}_{i}:=\mu_{k+1+i}\;\mbox{ for }\;i=1,\ldots,m
$$
and using the expressions for $\th^k_j$ in (\ref{theta}), we find from the above subgradients
\begin{eqnarray*}
\begin{array}{ll}
&(v_{j},w_{j})\in\partial f_{j}(\ox_{j},\oy_{j}),\;j=1,\ldots,k,\quad u^{k}_{i}\in\partial\varphi_{i}(\bar x_{k}),\;i=0,\ldots,m+r,\\\\
&\mbox{ and }\;u'^{k}_{i}\in\partial(-\varphi_{i})(\bar x_{k}),\;i=m+1,\ldots,m+r,
\end{array}
\end{eqnarray*}
for which we have the conditions
\begin{eqnarray*}
-x^*_{jj}=\lambda^{k}_0 h_{k}v_{j}+\psi_{j}-\psi_{j+1},\quad j=1,\ldots,k-1,
\end{eqnarray*}
\begin{eqnarray*}
-x^*_{k0}-x^*_{kk}=\lambda^{k}_0 h_{k}v_{k}+\psi_{k}+\sum^{m}_{i=0}\lambda^{k}_{0}u^{k}_{i}+\sum^{m+r}_{i=m+1}\mu_{k+1+i}u^{k}_{i}+\sum^{m+r}_{i=m+1}\mu_{k+1+r+i}u'^{k}_{i},
\end{eqnarray*}
\begin{eqnarray*}
-y^*_{jj}=\lambda^{k}_0 h_{k}w_{j}+\lambda^{k}_0 \theta^{k}_{j}-h_{k}\psi_{j},\quad j=1,\ldots,k.
\end{eqnarray*}
Next we introduce for each $k\in\N$ the adjoint discrete trajectories by
$$
p^{k}_{j-1}:=\psi^{k}_{j}\;\mbox{ for }\;j=1,\ldots, k\quad\mbox{and}
$$
$$
p^{k}_{k}:=-x^*_{k0}-\disp\sum^{m}_{i=0}\lambda^{k}_{i}u^{k}_{i}-\disp\sum^{m+r}_{i=m+1}\mu_{k+1+i}u^{k}_{i}-\disp\sum^{m+r}_{i=m+1}\mu_{k+1+r+i}u'^{k}_{i}.
$$
Then we get the relationships
$$
\frac{p^{k}_{j}-p^k_{j-1}}{h_k}=\frac{\psi^{k}_{j+1}-\psi^k_j}{h_k}=\lambda^{k}_{0}v_{j}+\frac{x^*_{jj}}{h_{k}},
$$
$$
p^{k}_{j-1}-\frac{\lambda^{k}_{0}\theta^{k}_{j}}{h_{k}}=\psi^{k}_{j}-\frac{\lambda^{k}_{0}\theta^{k}_{j}}{h_{k}}=\lambda^{k}_{0}w_{j}+\frac{y^*_{jj}}{h_{k}},
$$
which ensure the validity of the extended Euler-Lagrange inclusion of the theorem for each $j=1,\ldots,k$. Furthermore, it follows from
(\ref{phimr}), (\ref{phir}) and the complementary slackness conditions in (\ref{slack}) that we have
$$
\mu_{k+1+j}(\varphi_{j}(x^k_k)-L\eta_k)=0\;\mbox{ and }\;\mu_{k+1+r+j}(-\varphi_{j}(x^k_k)-L\eta_k)=0\;\mbox{ for }\;j=m+1,\ldots,m+r,
$$
which implies that either $\mu_{k+1+j}=0$ or $\mu_{k+1+r+j}=0$ must be equal to zero for all $j=m+1,\ldots,m+r$. Denoting finally
\begin{eqnarray}\label{lam}
\lambda^k_i:=\left\{\begin{array}{ll}
\mu_{k+1+i}&\mbox{if }\;\mu_{k+1+r+i}=0,\\
-\mu_{k+1+r+i}&\mbox{if }\;\mu_{k+1+i}=0
\end{array}\right.
\end{eqnarray}
for each $i=m+1,\ldots,m+r$, we get
\begin{eqnarray*}
\begin{array}{lll}
-p^{k}_{k}&=x^*_{k0}+\disp\sum^{m}_{i=0}\lambda^{k}_{i}u^{k}_{i}+\disp\sum^{m+r}_{i=m+1}\mu_{k+1+i}u^{k}_{i}+\disp\sum^{m+r}_{i=m+1}\mu_{k+1+r+i}u'^{k}_{i}\\
&\in N(\bar{x}^{k}_{k};\Omega_k)+\disp\sum^{m}_{i=0}\lambda^{k}_{i}\partial\varphi_{i}(\bar{x}^{k}_{k})+\disp\sum^{m+r}_{i=m+1}\mu_{k+1+i}\partial\varphi_{i}
(\bar{x}^{k}_{k})+\disp\sum^{m+r}_{i=m+1}\mu_{k+1+r+i}\partial(-\varphi_{i})(\bar{x}^{k}_{k})\\
&\subset\disp\sum^{m}_{i=0}\lambda^{k}_{i}\partial\varphi_{i}(\bar{x}^{k}_{k})+\disp\sum^{m+r}_{i=m+1}\lambda^{k}_{i}\partial^{0} \varphi_{i}(\bar{x}^{k}_{k})+N(\bar{x}^{k}_{k};\Omega_k).
\end{array}
\end{eqnarray*}
This justifies the transversality inclusion completes the proof of the theorem. $\h$\vspace*{0.05in}

The last result of this section specifies the nontriviality condition of Theorem~\ref{dnc} (meaning that all the dual elements therein, i.e., $\lm^k_i$ for $i=0,\ldots,m+r$ and $p^k_j$ for $j=0,\ldots,k$, are not equal to zero simultaneously) for the important class of multifunctions $F_j=F(\cdot,t_j)$ in the discrete inclusions (\ref{ie}) of the implicit Euler scheme satisfying the so-called Lipschitz-like (known also as Aubin's pseudo-Lipschitz) property around the optimal solution $\ox^k$ for $(P_k)$. Recall that a set-valued mapping $F\colon\R^n\tto\R^m$ is {\em Lipschitz-like} around $(\ox,\oy)\in\gph F$ if there exist neighborhoods $U$ of $\ox$ and $V$ of $\oy$ as well as a constant $\kappa\ge 0$ such that we have the inclusion
$$
F(u)\cap V\subset F(x)+\kappa|x-u|\B\;\mbox{ for all }\;x,u\in U.
$$
A crucial advantage of the nonconvex normal cone (\ref{mnc}) is the possibility to obtain in its terms a complete characterization of the Lipschitz-like property of arbitrary closed-graph multifunctions. To formulate this result, we recall coderivative notion for set-valued mappings generated by the normal cone (\ref{mnc}). Given $F\colon\R^n\tto\R^m$ and $(\ox,\oy)\in\gph F$, the {\em coderivative} of $F$ at $(\ox,\oy)$ is a set-valued mapping $D^*F(\ox,\oy)\colon\R^m\tto\R^n$ defined by
\begin{eqnarray}\label{cod}
D^*F(\ox,\oy)(v):=\Big\{u\in\R^n\Big|\;(u,-v)\in N((\ox,\oy);\gph F)\Big\}\;\mbox{ for all }\;v\in\R^m.
\end{eqnarray}
When $F$ is single-valued and smooth around $\ox$ (then we drop $\oy=F(\ox))$, its coderivative reduces to the adjoint/transpose Jacobian
$$
D^*F(\ox)(v)=\{\nabla F(\ox)^Tv\},\quad v\in\R^m.
$$
In the general nonsmooth and/or set-valued case, the coderivative (\ref{cod}) is a positive homogeneous multifunction, which enjoys comprehensive calculus rules based on the variational and extremal principle of variational analysis; see \cite{bs1,rw}.

The results we need in what follows in known as the {\em coderivative/Mordukhovich criterion} (see \cite[Theorem~5.7]{bs93} and \cite[Theorem~9.40]{rw} with the references therein): If $F$ is close-graph around $(\ox,\oy)$, then it is Lipschitz-like around this point if and only if
\begin{eqnarray}\label{cc}
D^*F(\ox,\oy)(0)=\{0\}.
\end{eqnarray}

Now we are ready to derive the aforementioned consequence of Theorem~\ref{dnc}.

\begin{Corollary}[enhanced nontriviality condition]\label{cp1} In addition to the assumptions of Theorem~{\rm(\ref{dnc})}, suppose that for each $j=1,\ldots,k$, the multifunction $F_{j}$ is Lipschitz-like around the optimal point $(\ox^{k}_{j},(\ox^k_j-\ox^k_{j-1})/h_k)$. Then all the necessary optimality conditions of this theorem hold at $\ox^k$ with the enhanced nontriviality
\begin{equation}\label{p1}
\disp\sum_{i=0}^{m+r}|\lambda^{k}_i|+|p^k_0|=1\;\mbox{ for all }\;k\in\N.
\end{equation}
\end{Corollary}
{\bf Proof.} If $\lm^k_0=0$, then it follows from the Euler-Lagrange inclusion of the theorem that
$$
\disp\Big(\frac{p^k_{j}-p^k_{j-1}}{h_k},p^{k}_{j-1}\Big)\in N\Big(\Big(\ox^{k}_{j},\frac{\ox^{k}_{j}-\ox^k_{j-1}}{h_k}\Big);\gph F_j\Big)
$$
for all $j=1,\ldots,k$, which tells us by the coderivative definition (\ref{cod}) that
\begin{eqnarray*}
\disp\frac{p^{k}_{j}-p^k_{j-1}}{h_k}\in D^*F_j\Big(\ox^{k}_{j},\frac{\ox^{k}_{j}-\ox^k_{j-1}}{h_k}\Big)(-p^k_{j-1}),\quad j=1,\ldots,k.
\end{eqnarray*}
Employing finally the coderivative criterion (\ref{cc}) with taking into account the transversality condition of the theorem as well as the normalization of $(\lm_0,\ldots,\lm_{m+r},p^k_0)$ without changing other conditions, we arrive at (\ref{p1}) and thus completes the proof. $\h$

\section{Concluding Remarks}\sce

This paper develops a constructive approach to investigate the generalized Bolza problem of optimizing constrained differential inclusions satisfying the relaxed one-sided Lipschitzian condition by using the implicit Euler scheme of discrete approximations. In this way we not only justify the well-posedness of the suggested discrete approximation procedures in the sense of either the uniform or $W^{1,2}$-convergence of discrete optimal solutions to a given local (strong or intermediate) minimizer of the original nonsmooth Bolza problem, but also derive necessary optimality conditions to solve each problem of the implicit Euler discrete approximations. As mentioned in the introductory Section~1, the results obtained are new even in the case of the implicit Euler scheme for unconstrained differential inclusions satisfying the classical Lipschitz condition.

A natural question arises about the possibility to derive {\em necessary optimality conditions} for the given intermediate or strong local minimizer of the original problem $(P)$ for favorable classes of {\em ROSL differential inclusions} by passing to the limit from those obtained for $(P_k)$ and $(\tilde P_k)$, respectively, as $k\to\infty$. It can surely be done in the case when, in the setting of Theorem~\ref{dsa}, the velocity function $F$ is Lipschitzian around the local minimizer under consideration; cf.\ \cite[Theorem~6.1]{bs95} and \cite[Theorem~6.22]{bs2}, where the case of the explicit Euler scheme was investigated. Note to this end that the ROSL and Lipschitz-like properties of $F$ used in Corollary~\ref{cp1} are generally independent (even for bounded mappings), and they both are implied by the classical Lipschitz condition.

On the other hand, the method of discrete approximations has been successfully employed in \cite{chhm1} to derive necessary optimality conditions for the Bolza problem governed by a dissipative (hence ROSL while unbounded and heavily non-Lipschitzian) differential inclusion that arises in optimal control of Moreau's {\em sweeping process} with mechanical applications. The procedure in \cite{chhm1} exploits some specific features of the controlled sweeping process over convex polyhedral sets, and thus a principal issue of the our further research is about the possibility to extend these results to more general ROSL differential inclusions.


\small

\begin{thebibliography}{10}

\bibitem{a} Z. Artstein, First-order approximations for differential inclusions, {\em Set-Valued Anal.} {\bf 2} (1994), 7–-18.

\bibitem{abm} H. Attouch, G. Buttazzo and G. Michaille, {\em Variational Analysis in Sobolev and BV Spaces}, SIAM Publications, Philadelphia,
PA, 2005.

\bibitem{br} W. J. Beyn and J. Rieger, The implicit Euler scheme for one-sided Lipschitz differential inclusions, {\em Discr. Cont. Dyn. Syst.}
{\bf 14} (2009), 409--428.

\bibitem{chhm} G. Colombo, R. Henrion, N. D. Hoang and B. S. Mordukhovich, Discrete aproximations of a control sweeping process, {\em Set-Valued
Var. Anal.} {\bf 23} (2015), 69--86.

\bibitem{chhm1} G. Colombo, R. Henrion, N. D. Hoang and B. S. Mordukhovich,  Optimal control of the sweeping process over polyhedral controlled sets,
preprint (2015).

\bibitem{td91} T. Donchev, Functional differential inclusions with monotone right-hand side, {\em Nonlinear Anal.} {\bf 16} (1991), 543--552.

\bibitem{td02} T. Donchev, Properties of one-sided Lipschitz mutivalued maps, {\em Nonlinear Anal.} {\bf 49} (2002), 13--20.

\bibitem{tdf98} T. Donchev and E. Farkhi, Stability and Euler approximations of one-sided Lipschitz convex differential inclusions,
{\em SIAM J. Control Optim.} {\bf 36} (1998), 780--796.

\bibitem{tdeb} T. Donchev, E. Farkhi and B. S. Mordukhovich, Discrete approximations, relaxation, and optimization of one-sided
Lipschitzian differential inclusions in Hilbert spaces, {\em J. Diff. Eqns.} {\bf 243} (2007), 301--328.

\bibitem{dfr} T. Donchev, E. Farkhi and S. Reich, Fixed set iterations for relaxed Lipschitz multimaps, {\em Nonlinear Anal.} {\bf 53} (2003), 997--1015.

\bibitem{dl} A. L. Dontchev and F. Lempio, Difference methods for differential inclusions: a survey, {\em SIAM Rev.} {\bf 34} (1992), 263--294.

\bibitem{lv} F. Lempio and V. Veliov, Discrete approximations of differential inclusions, {\em Bayreuth. Math. Schr.} {\bf 54} (1998), 149–-232.

\bibitem{bs76} B. S. Mordukhovich, Maximum principle in problems of time optimal controls with nonsmooth constraints, {\em J. Appl. Math. Mech.}
{\bf 40} (1976), 960--969.

\bibitem{bs95} B. S. Mordukhovich, Discrete approximations and refined Euler-Lagrange conditions for differential inclusions,
{\em SIAM J. Control Optim.} {\bf 33} (1995), 882--915.

\bibitem{bs93} B. S. Mordukhovich, Complete characterization of openness, metric regularity, and Lipschitzian properties
of multifunctions, {\em Trans. Amer. Math. Soc.} {\bf 340} (1993), 1--35.

\bibitem{bs1} B. S. Mordukhovich, {\em Variational Analysis and Generalized Differentiation, I: Basic Theory},  Springer, Berlin, 2006.

\bibitem{bs2} B. S. Mordukhovich, {\em Variational Analysis and Generalized Differentiation, II: Applications}, Springer,
Berlin, 2006.

\bibitem{rw} R. T. Rockafellar and R. J-B. Wets, {\em Variational Analysis}, Springer, Berlin, 1998.

\bibitem{s} G. V. Smirnov, {\em Introduction to the Theory of Differential Inclusions}, AMS Publications, Providence, RI, 2001.

\bibitem{v00} R. B. Vinter, {\em Optimal Control}, Birkh\"auser, Boston, MA, 2000.

\end{thebibliography}
\end{document}